\documentclass[12pt]{amsart}
\newtheorem{thm}{Theorem}
\newtheorem{prop}[thm]{Proposition}
\newtheorem{lem}[thm]{Lemma}

\theoremstyle{remark}
\newtheorem{rem}[thm]{Remark}
\newtheorem{ex}[thm]{Example}

\theoremstyle{definition}
\newtheorem{defn}[thm]{Definition}

\newcommand{\thmref}[1]{Theorem~\ref{#1}}

\newcommand{\propref}[1]{Proposition~\ref{#1}}

\newcommand{\C}{\mathbb{ C}}

\newcommand{\lra}{\longrightarrow}

\newcommand{\GG}{\mathcal{ G}}
\newcommand{\FF}{\mathcal{ F}}

\newcommand{\R}{\mathbb{ R}}

\newcommand{\Id}{\operatorname{Id}}

\newcommand{\Z}{\mathbb{ Z}}
\newcommand{\cp}{{\C P^2}}

\newcommand{\x}{\times}

\newcommand{\al}{\alpha}
\newcommand{\bb}{\beta}
\newcommand{\cc}{\gamma}
\newcommand{\ww}{\wedge}
\newcommand{\lk}{\mathrm{lk}}

\title{Linking numbers of measured foliations}
\author{D.~Kotschick}
\address{Mathematisches~Institut,~Ludwig-Maximilians-Universit\"at
M\"unchen, Theresienstr.~39, 80333 M\"unchen, Germany}
\email{dieter@member.ams.org}
\author{T.~Vogel}
\address{Mathematisches~Institut,~Ludwig-Maximilians-Universit\"at
M\"unchen, Theresienstr.~39, 80333 M\"unchen, Germany}
\email{thomas.vogel@mathematik.uni-muenchen.de}

\date{\today; MSC 2000: 57R30, 37C10, 37C40, 37C85.}

\thanks{The first author gratefully acknowledges support from the
{\it VolkswagenStiftung} and the hospitality of Harvard University.
The second author is supported by the {\it DFG Graduiertenkolleg
``Mathematik im Bereich ihrer Wechselwirkung mit der Physik''.}
}

\begin{document}

\begin{abstract}
We generalise the average asymptotic linking number of a pair of
divergence-free vector fields on homology three-spheres~\cite{A,AK,V}
by considering the linking of a divergence-free vector field on a
manifold of arbitrary dimension with a codimension two foliation
endowed with an invariant transverse measure. We prove that the
average asymptotic linking number is given by an integral of Hopf type.
Considering appropriate vector fields and measured foliations, we
obtain an ergodic interpretation of the Godbillon-Vey invariant of a
family of codimension one foliations discussed in~\cite{K}.
\end{abstract}

\maketitle

\section{Introduction}\label{s:intro}

Generalising the differential form description of the classical Hopf
invariant, various authors in hydrodynamics introduced a Hopf invariant
$H(X,Y)\in\R$ for pairs $(X,Y)$ of divergence-free vector fields on a
three-dimensional homology sphere $M$. Arnold~\cite{A} interpreted this
Hopf invariant as the average over $M\times M$ of the long-time asymptotic 
linking numbers of pairs of flow lines for $X$ and $Y$. In order to make 
this interpretation precise, one has to find a coherent way of closing 
long pieces of flow lines to form loops so that one can evaluate the 
corresponding linking numbers and compute their long-time asymptotics. 
Arnold's construction of such a system of short path works for suitably 
generic vector fields, cf.~\cite{AK}, but there are technical problems 
if one allows vector fields with degenerate zeroes. A slightly different 
construction covering all divergence-free vector fields was given by the 
second author in~\cite{V}. 

A variation of Arnold's construction formulated for nonsingular flows 
on $S^{3}$ with invariant measures is contained in~\cite{CI}.
Higher-dimensional generalisations are somewhat
elusive, because it is not at all clear how to parameterise and ``close''
higher-dimensional open submanifolds in order to define linking
numbers. In this paper we define the average asymptotic linking number 
of a divergence-free vector field and a measured codimension two foliation. 
We close up pieces of flow lines of the vector field as in~\cite{V}, but 
we do not parameterise or close the leaves of the foliation. Intuitively, 
parameterisation and ``closing'' is done by the invariant measure. Thinking 
of a divergence-free flow as a (singular) one-dimensional foliation with 
a holonomy-invariant transverse measure, we recover the construction 
of~\cite{A,V}. A different generalisation of Arnold's construction
has been proposed by Rivi\`ere~\cite{R}.

In section~\ref{s:sub} we shall consider the linking of the flow lines
of a divergence-free vector field with a null-homologous closed oriented
submanifold $N$ of codimension two. On three-manifolds, this situation 
was already considered by Arnold~\cite{A}, but our approach is different. 
For manifolds with vanishing first Betti number we prove the existence of 
suitable systems of short paths which one can use to close the flow lines 
of the vector field. The resulting average asymptotic linking number is 
given by a Hopf-type integral. It is then clear how to generalise further 
and replace the closed submanifold $N$ by a measured foliation $(\FF,\nu)$ 
of codimension two, whose Ruelle--Sullivan cycle is null-homologous. 
This is carried out in section~\ref{s:measures}. The discussion there 
is motivated in part by the work of Arnold, Khesin~\cite{Kh}, Novikov 
and others on higher-dimensional generalisations of Arnold's construction 
described in~\cite{AK} (Chapter III, 7.B). 

In section~\ref{s:GV} we apply our construction to give an
interpretation of the four-dimensional Godbillon-Vey invariant
discussed in~\cite{K} as an average asymptotic linking number of a
vector field and a measured codimension two foliation. For this it is 
important that the constructions of section~\ref{s:measures} work
for singular foliations.

\medskip
\noindent
{\sl Acknowledgements:} We are grateful to B.~Khesin for useful 
conversations, for Remark~\ref{r:Khesin}, and for sending us a copy 
of~\cite{Kh}; and to S.~Hurder for prodding us to be more careful in our 
treatment of singular foliations, and for pointing us to the work of Moussu.

\section{Preliminaries}\label{s:prelim}

Let $M^{n}$ be a smooth closed oriented $n$-manifold.
Consider smooth closed oriented submanifolds $N_{1}^{k}$ and
$N_{2}^{l}$ of $M$, with $k+l=n-1$. If $N_{1}$ and $N_{2}$ are
null-homologous over $\R$, they have a well-defined linking number
defined as follows. Let $S_{1}$ be a real oriented $(k+1)$-chain with
$\partial S_{1}=N_{1}$, and define $\lk(N_{1},N_{2})$ to be the
algebraic intersection number of $S_{1}$ and $N_{2}$, which by
assumption have complementary dimensions. This linking number
is always rational. If $N_{1}$ and $N_{2}$ are
null-homologous over $\Z$, then their linking number is integral.
In any case, we have $\lk(N_{1},N_{2})=(-1)^{(k+1)(l+1)}\lk(N_{2},N_{1})$.

There is a classical expression for these linking numbers as integrals
of certain differential forms. If $W_{i}$ is a tubular neighbourhood
of $N_i$, we can choose a closed form $\eta_{i}$ on $W_{i}$ representing
the compactly supported Poincar\'e dual of $N_{i}$ in $W_{i}$, and then
extend this form by zero to all of $M$. The assumption that $[N_{1}]$
vanishes in real homology means that the extension of $\eta_{1}$ to
$M$, also denoted $\eta_{1}$, is exact. Let $\alpha_{1}$
be a primitive and define
\begin{equation}\label{eq:Hopf}
H(N_{1},N_{2})=\int_{M}\alpha_{1}\ww\eta_{2} \ .
\end{equation}
If we chose a different primitive $\alpha_{1}'$ for $\eta_{1}$, then
$$
\int_{M}\alpha_{1}'\ww\eta_{2}-\int_{M}\alpha_{1}\ww\eta_{2}=
\int_{M}(\alpha_{1}'-\alpha_{1})\ww\eta_{2} \ .
$$
This vanishes because $\alpha_{1}'-\alpha$ is closed and $\eta_{2}$
is exact as $N_{2}$ is also assumed to be null-homologous. If we make
a different choice, $\eta_{2}'$, for $\eta_{2}$, then
$\eta_{2}'-\eta_{2}=d\beta$, for a $\beta$ with support in $W_{2}$.
Then
$$
\int_{M}\alpha_{1}\ww\eta_{2}'-\int_{M}\alpha_{1}\ww\eta_{2}=
\int_{M}\alpha_{1}\ww d\beta=\pm\int_{M}d\alpha_{1}\ww\beta=
\pm\int_{M}\eta_{1}\ww\beta \ .
$$
Now the right hand side vanishes because $\eta_{1}$ and $\beta$ have
support in $W_{1}$ and $W_{2}$ respectively, which can be chosen to
be disjoint.

Thus we have shown that $H(N_{1},N_{2})$ is independent of choices,
except possibly for the choice of $\eta_{1}$. However, we can interchange 
the roles of $\eta_{1}$ and $\eta_{2}$ because we have $H(N_{1},N_{2})=
(-1)^{(k+1)(l+1)}H(N_{2},N_{1})$. Thus $H(N_{1},N_{2})$ is
independent of all choices and only depends on the submanifolds
$N_{1}$ and $N_{2}$, as suggested by the notation. By making a 
convenient choice for $\eta_{2}$ as in~\cite{BT}, one sees that
\begin{equation}\label{eq:link}
\int_{M}\alpha_{1}\ww\eta_{2} = \int_{N_{2}}\alpha_{1} 
\end{equation}
and that these integrals equal the linking number $\lk(N_{1},N_{2})$. 
An instance of this equality is the differential form interpretation 
of the Hopf invariant, and we shall refer to either 
$\int_{M}\alpha_{1}\ww\eta_{2}$ or $\int_{N_{2}}\alpha_{1}$ as a 
Hopf-type integral.

The linking number of two disjoint closed oriented submanifolds
can be expressed through so-called linking forms (see~\cite{deR,V}).
A {\em linking form} $L$ on an $n$-manifold $M$ is a double form on
$M\times M$ with the property that whenever the linking number of two
oriented closed submanifolds $N_1$ and $N_2$ is well-defined, we have
$$
\lk(N_1,N_2) = \int_{N_1} \int_{N_2} L \ .
$$
A linking form can be constructed as follows. We choose a Riemannian
metric on $M$, and define $H$ to be the projection operator mapping a
differential form to its harmonic part. Then, for every degree $i$,
the Green's operator
$$
G \colon \Omega^{i}(M)\longrightarrow ({\mathcal H}^{i})^{\perp}
$$ 
is defined to map a differential form $\alpha$ of degree $i$ to the 
unique $i$-form $\omega$ that is perpendicular to all harmonic forms 
and solves the equation
$$
\Delta\omega = \alpha - H(\alpha) \ .
$$
In other words, $G$ is characterised by the properties $HG=0$ and
$\Delta G=\Id -H$. It can be written in the form
$$
G(\alpha)(x)=\int_{y\in M}\alpha (y)\ww *_{y}g(x,y)
$$
for a suitable double form $g(x,y)$ on $M\times M$ which is smooth
away from the diagonal and has a pole of order $n-2$ along the 
diagonal, cf.~\cite{deR} \S 31.
Here $*_y$ denotes the Hodge star operator with respect to the second
factor of $M\x M$. 

We denote by $(-1)^\epsilon$ the linear operator on double forms 
acting on decomposable double forms $\omega(x)\cdot\eta(y)$ with 
$\eta\in\Omega^s(M)$ by multiplication with 
$(-1)^{(n-s)s}$. In later sections we will be concerned only with the 
case $s=n-1$. (This is due to the fact that the Poincar\'e duals of 
one-dimensional submanifolds of $M$ have degree $n-1$.)
We define the double form $L$ on $M\x M$ by
\begin{equation}\label{e:defL}
L(x,y) = (-1)^\epsilon *_yd_yg(x,y) \ .
\end{equation}

\begin{prop}\label{p:fundL}
The double form $L(x,y)$ is a linking form. Denoting by $r$ the 
Riemannian distance function, $L$ has a singularity of order 
$(r(x,y))^{1-n}$ along the diagonal in $M\x M$ and is smooth 
elsewhere. It has the following additional property: for every 
$i$-form $\al$ there exists an $(i-1)$-form $h$ such that
\begin{equation}\label{e:fundL}
\int_{y\in M} L(x,y)\ww d\al(y) = \al(x) -H(\al)(x)+dh(x) \ . 
\end{equation}
\end{prop}
\begin{proof}
By the definition of $G$, and because $G$ commutes with $\Delta$, we have
$$
G(d^{*}d\al)=\al-H(\al)-G(dd^{*}\al) \ .
$$
As $G$ commutes with $d$ we can set $h=-G(d^{*}\al)$ to obtain
\begin{equation*}
\int_{y\in M} d\al(y)\ww *_yd_yg(x,y) = G(d^*d\al)= 
\al(x) -H(\al)(x)+ dh(x) \ . 
\end{equation*}
If we change the order of the factors in the integrand, we have to multiply 
by $(-1)^\epsilon$ and we obtain~\eqref{e:fundL}. 

That $L$ is a linking form can be shown as in the proof of Theorem 3
in~\cite{V}. We briefly indicate the argument. Choosing $W_{i}$,
$\eta_{i}$ and $\alpha_{i}$ as above, and using~\eqref{e:fundL}, we have
\begin{align*}
lk(N_1 &,N_2)=\int_M \alpha_1\ww\eta_2 \\
           =&\int_{x\in M}\left(H(\alpha_{1})(x)-dh(x)+
       \int_{y\in M}L(x,y) \ww \eta_1(y) \right)\ww\eta_2(x) \ .
\end{align*}
The integrals
$$
\int_{M}H(\alpha_{1})\ww\eta_{2} \quad \mathrm{and} \quad
\int_{M}dh\ww\eta_{2}
$$
vanish by Stokes's theorem because $\eta_{2}$ is exact. As $\eta_{i}$
has support in $W_{i}$, we obtain
\begin{align*}
lk(N_1,N_2)=&\int_{x\in W_{2}}\left(\int_{y\in W_{1}}
L(x,y)\ww\eta_1(y)\right)\ww\eta_2(x) \ .
\end{align*}
Now the $W_{i}$ can be taken to be disjoint and making a good choice for 
the $\eta_{i}$ the above integral reduces to
\begin{equation*}
\int_{x\in N_2} \int_{y\in N_1} L(x,y) \ . 
\end{equation*}
The claim about the order of the singularity along the diagonal 
follows from what we said above about the singularity of $g$.
\end{proof}

We shall also need to use the mean ergodic theorem as in~\cite{V}.
\begin{thm}[\cite{dus}]\label{t:ergodic}
Let $f$ be an $L^1$-function on the compact manifold $M$. Let $\phi_t$ be a
differentiable flow on $M$ preserving a given volume form $\mu$.
\begin{enumerate}
\item The limit
$$
\tilde{f}(x)=\lim_{t\to \infty} \frac{1}{t}\int_{s=0}^t f(\phi_s(x)) ds
$$
exists in the $L^1$-sense and is an integrable function.
\item The integral of $\tilde{f}$ satisfies
$$
\int_M \tilde{f} \mu = \int_M f \mu\ .
$$
\end{enumerate}
\end{thm}

\section{Linking numbers between vector fields and submanifolds}\label{s:sub}

Let $M^{n}$ be a smooth closed oriented $n$-manifold, with $n\geq 3$.
We fix once and for all a Riemannian metric and the corresponding linking
form $L$ as in Proposition~\ref{p:fundL}.

Instead of considering the linking between two closed oriented submanifolds
$N_{1}$ and $N_{2}$, we shall replace $N_{1}$ by loops formed by closing
up the flow lines of a divergence-free vector field $X$, and consider how
these loops link with an oriented submanifold $N^{n-2}$ playing the role
of $N_{2}$ above. As before, we assume that $N\subset M$ is
null-homologous over $\R$. To ensure that all the loops are
null-homologous, we assume that the first Betti number of $M$
vanishes.

Let $\mu$ be a volume form on $M$, and $X$ a vector field that is
divergence-free with respect to $\mu$, i.~e.~such that $L_{X}\mu=0$.
Then the $(n-1)$-form $\eta=i_{X}\mu$ is closed, and our assumption
that the first Betti number of $M$ vanishes implies, via Poincar\'e
duality, that $\eta$ must be exact. Let $\alpha$ be a primitive.
Then, generalising~\eqref{eq:link}, we can define a Hopf-type integral
for $X$ and $N$ by setting
$$
H(X,N)=\int_{N}\alpha \ .
$$
If we choose a different primitive $\alpha'$ for $\eta=i_{X}\mu$, then
$$
\int_{N}\alpha'-\int_{N}\alpha = \int_{N}\alpha'-\alpha
$$
vanishes because $\alpha'-\alpha$ is closed and $N$ is assumed to be
homologous to zero. Thus $H(X,N)$ is well-defined.

We want to interpret this integral as an average of asymptotic
linking numbers of flow lines of $X$ with $N$. To do so, we fix
once and for all a ``system of short paths'' connecting any pair of
points $p,q\in M$. The flow of $X$ will be denoted by $\phi_t$.

\begin{defn}\label{d:spsN}
A set $\Sigma$ of piecewise differentiable paths in $M$ is {\em a system
of short paths} if it has the following properties:
\begin{enumerate}
\item For any two points $p,q\in M$ there is exactly one path $\sigma(p,q)
\in\Sigma$ starting at $p$ and ending at $q$.
\item The paths depend continuously on their endpoints almost everywhere.
\item The limit
\begin{equation}\label{spscondN}
\lim_{t\to\infty} \frac{1}{t} \int_{y\in\sigma(\phi_t(x),x)}\int_{p\in N}
L(p,y) = 0
\end{equation}
exists in the $L^1$-sense.
\end{enumerate}
\end{defn}

\begin{thm}\label{t:spsN}
A system of short paths exists. It can be chosen independently of the
vector field $X$.
\end{thm}
\begin{proof}
Let $C(y)=\int_{p\in N} L(p,y)$. This is a well-defined $L^{1}$-form on $M$.
Since $M$ is compact it can be covered by a finite number of geodesic balls 
$U_{j}$, $j=1,\ldots,r$. In each geodesic ball we fix a basepoint 
$u_{j}\in U_{j}$ such that the following conditions are satisfied.
\begin{enumerate}
\item For every pair $k,j$ there is a path $\gamma_{kj}$ parametrised 
by $[0,1]$ joining $u_{k}$ and $u_{j}$ such that the integral
\begin{equation}\label{e:intijalt}
\int_{y\in\gamma_{kj}} |C(y)| \ := \  
\int_0^1\left|i_{\frac{\partial\gamma_{kj}}{\partial s}}C(\gamma_{kj}(s))
\right| ds
\end{equation}
is finite.
\item Let $\sigma(x,u_{j})$ denote the unique geodesic in $U_{j}$ between 
$x\in U_{j}$ and $u_{j}$. For all $j$ the integral
\begin{equation}
\int_{x\in U_{j}} \int_{y\in\sigma(x,u_{j})} |C(y)|\mu(x)
\end{equation} 
is finite.
\end{enumerate}
The second condition is satisfied if all the $u_{j}$ are outside of $N$. 
Since $N$ has codimension two, both conditions hold for a generic choice 
of the $u_{j}$.

For all $x\in M$ fix a number $n(x)$ such that $x\in U_{n(x)}$ and $n(x)$ is 
locally constant on a dense open subset of $M$. Let $p,q\in M$. We define a 
piecewise differentiable path $\sigma(p,q)$ joining $p$ and $q$ as follows. 
The first segment of $\sigma(p,q)$ is the unique geodesic in $U_{n(p)}$ between 
$p$ and $u_{n(p)}$, the second segment is $\gamma_{n(p),n(q)}$ and the third is 
the unique geodesic in $U_{n(q)}$ starting at $u_{n(q)}$ and ending at $q$. 
Note that for $x$ with $n(x)=j$, the short path between $x$ and $u_j$ is the 
unique geodesic in $U_j$ joining $x$ and $u_j$, this justifies the notation 
$\sigma$ for both objects.

We define $\Sigma$ to be the set of paths obtained this way. For each 
$p,q\in M$ there is a unique path with starting point $p$ and end point $q$. 
By construction, the paths depend in a continuous way on their starting and 
end points on an open dense subset of $M\x M$. 

Dividing the paths into their differentiable 
pieces we obtain
\begin{align*}
\int_{x\in M} \left|\int_{y\in\sigma(\phi_t(x),x)} C(y)\right|\mu(x) & 
\le \int_{x\in M} \int_{y\in\sigma(\phi_t(x), u_{n(\phi_t(x))})} |C(y)| \mu(x)\\
& \quad + \int_{x\in M} \int_{y\in \gamma_{n(\phi_t(x)),n(x)}} |C(y)|
\mu(x) \\ 
& \quad + \int_{x\in M} \int_{y\in \sigma(u_{n(x)},x)} |C(y)| \mu(x)~.
\end{align*}
The first and the third summand on the right hand side are in fact equal. 
To see this, apply the volume-preserving transformation $x\mapsto \phi_{-t}(x)$ 
to the first summand. In particular these two summands do not depend on $t$. 
The second summand is bounded above by 
\begin{equation*}
\max_{i,j}\left(\int_{y\in\gamma_{ij}} |C(y)|\right) \ \int_{x\in M} \mu(x) 
\end{equation*}
and the third summand is bounded above by
\begin{equation*}
\sum_{i=1}^r \int_{x\in U_i} \int_{y\in\sigma(u_i,x)} |C(y)| \mu(x)~.
\end{equation*}
Thus $\Sigma$ meets the third condition of Definition~\ref{d:spsN}.  
\end{proof}

If $x$ is a point in $M$ and $t\in\R$, denote by $\phi(x,t)$ the
flow line of $X$ generated by $x$ in the time-interval $[0,t]$. We
shall denote by $\gamma(x,t)$ the closed loop obtained by connecting
the endpoints $x$ and $\phi_{t}(x)$ of $\phi(x,t)$ by the path
$\sigma(\phi_{t}(x),x)$.

The following is true for dimension reasons:
\begin{lem}
Let $t\in\R$ be fixed. Then for $\mu$-almost all $x\in M$, the piecewise
differentiable curve $\gamma(x,t)$ is embedded in $M\setminus N$.
\end{lem}
Given this Lemma and our assumptions that the first Betti
number of $M$ vanishes and that $N$ is null-homologous over $\R$,
we can define the linking number $\lk(\gamma(x,t),N)$ for almost
all $x$ and $t$. We then have:
\begin{prop}\label{p:asymptotic}
The limit 
$$
\lk(x,N)=\lim_{t\to\infty}\frac{1}{t}\lk(\gamma(x,t),N)
$$
exists in the $L^{1}$-sense. It is an integrable function on $M$
which does not depend on the chosen system of short paths.
\end{prop}
\begin{proof}
By Proposition~\ref{p:fundL} we have
\begin{equation*}
\lim_{t\to\infty}\frac{1}{t}\lk(\gamma(x,t),N)
= \lim_{t\to\infty}\frac{1}{t}\int_{y\in\gamma(x,t)}\int_{p\in N} L(p,y) \ .
\end{equation*}
Using the third property of the system of short paths, we find that
the right hand side equals
\begin{equation*}
    \lim_{t\to\infty}\frac{1}{t}\int_{s=0}^t\int_{p\in N}
         i_XL(p,\phi_sx) ds\ .
\end{equation*}
Since $L(p,\cdot)$ is an integrable form on $M$ for every $p\in M$,
we can apply the mean ergodic theorem, Theorem~\ref{t:ergodic}. Hence 
the limit exists. It is clearly independent of the system of short paths.
\end{proof}

Using this, we can finally define the average asymptotic linking
number of the vector field $X$ with the submanifold $N$, by setting
$$
\lk(X,N)=\int_{M}\lk(x,N) \mu \ .
$$
The analog of the theorem proved in~\cite{A,V} for vector fields on
$3$-manifolds is:
\begin{thm}\label{t:sub}
Let $M$ be a closed oriented $n$-manifold with $b_1(M)=0$.
Let $X$ be a divergence-free vector field on $M$, and $N\subset M$ a
closed oriented submanifold of codimension $2$ which is null-homologous
over $\R$. Then the average asymptotic linking number of the orbits of
$X$ with $N$ exists and equals a Hopf-type integral:
$$
\lk(X,N)=H(X,N) \ .
$$
\end{thm}
\begin{proof}
We have seen in the proof of \propref{p:asymptotic} that
\begin{align*}
\lk(x,N) = & \lim_{t\to\infty}\frac{1}{t}\int_{s=0}^t\int_{p\in N}
         i_XL(p,\phi_sx) ds\ .
\end{align*}
By the second part of the mean ergodic theorem, Theorem~\ref{t:ergodic}, 
we find
\begin{align*}
\lk(X,N) = & \int_{x\in M}
\left(\lim_{t\to\infty}\frac{1}{t}\int_{s=0}^t\int_{p\in N}
         i_XL(p,\phi_sx) ds \right) \mu(x) \\
     = & \int_{x\in M} \int_{p\in N} i_{X(x)}L(p,x)\ww \mu(x) \\
     = & \int_{x\in M} \int_{p\in N} L(p,x)\ww i_{X(x)}\mu(x)\ ,
\end{align*}
where $i_{X}\mu = d\alpha$. We now apply Proposition~\ref{p:fundL} to the last 
integral. The integral of the harmonic term in Proposition~\ref{p:fundL} 
over $N$ is zero because $N$ is null-homologous, and by Stokes's theorem the 
exact term plays no role. Hence we obtain the desired result
\begin{equation*}
\lk(X,N) = \int_N \alpha =H(X,N) \ .
\end{equation*}
\end{proof}

To end this section we consider some geometric examples.

\begin{ex} 
Let $(M,\omega)$ be a symplectic manifold of dimension $2n$. 
Then $\omega^n$ is a volume form on $M$. Consider a smooth function $H$ 
on $M$ and the corresponding Hamiltonian vector field $X_{H}$. This is the 
unique vector field on $M$ satisfying $dH=i_{X_{H}}\omega$. Then $X_{H}$ 
is divergence-free with respect to $\omega^n$ since
$$
L_{X_{H}}\omega^{n}=di_{X_{H}}\omega^n = d(ndH\ww\omega^{n-1}) = 0 \ .
$$
Moreover, the $(2n-2)$--form $nH\omega^{n-1}$ satisfies 
$$
d(nH\omega^{n-1}) = ndH\ww\omega^{n-1} = n(i_{X_{H}}\omega)\ww\omega^{n-1} 
= i_{X_{H}}(\omega^{n}) \ .
$$
Let $N$ be a null-homologous submanifold of codimension $2$ in $M$. 
If $b_{1}(M)=0$, the linking number of $X$ with $N$ is given by
\begin{equation}\label{e:symplbsp}
\mathrm{lk}(X,N)=n\int_N H\omega^{n-1} \ .
\end{equation}
For a fixed Hamiltonian vector field the function generating it
is well-defined only up to the addition of locally constant functions. 
By Stokes's theorem this ambiguity in the choice of $H$ does not change 
the value of the integral in~\eqref{e:symplbsp} since $N$ is homologous 
to zero.
\end{ex}

\begin{ex}
Consider a manifold of dimension $2n+1$ with a contact form 
$\al\in\Omega^1(M)$. This means that $\al\ww(d\al)^n$ is a volume form 
on $M$. There is a unique vector field $X$, called the Reeb vector 
field, with the properties $\al(X)=1$ and $i_Xd\al=0$. Because
$$
L_{X}(\al\ww(d\al)^n) =di_X(\al\ww(d\al)^n) = d (d\al)^n = 0~,
$$ 
the Reeb vector field is divergence--free  with respect to the volume form 
$\al\ww(d\al)^n$. The $(2n-1)$--form $\al\ww(d\al)^{n-1}$ is a primitive of 
$i_X(\al\ww(d\al)^n) = (d\al)^n$. Thus, if $b_{1}(M)=0$, the linking number 
of the Reeb vector field $X$ with a null-homologous submanifold $N$ of 
codimension two is
$$
\mathrm{lk}(X,N) = \int_N \al\ww(d\al)^{n-1} \ .
$$
In particular it is nonzero if $N$ is a contact submanifold.
\end{ex}

\begin{ex}
Consider $M = S^3\x S^3\subset \C^2\x\C^2$ and the map 
\begin{alignat*}{1}
f \colon S^3\x S^3 & \lra \C \\
     (z,w) & \longmapsto \langle z,w \rangle = z_0\bar{w}_0 + 
     z_1\bar{w}_1 \ .
\end{alignat*}
As $0$ is a regular value of $f$ the preimage $N=f^{-1}(0)$ is a smooth 
submanifold of $M$ of dimension four. The map 
\begin{alignat*}{1}
S^3\x S^1 &\lra N \\
((z_0,z_1),\lambda) &\longmapsto ((z_0,z_1),(\lambda\bar{z}_1, -\lambda\bar{z}_0))
\end{alignat*}
is a diffeomorphism. Since $H_4(M)=0$, $N$ is null-homologous. We now 
construct a submanifold in $M$ whose boundary is $N$. Consider the set 
$[0,\infty)\subset\C$ of nonnegative real numbers. Let $(z,w)$ be a point 
in $M$ that is mapped by $f$ to a positive real number $c$. Then
$$
\left.\frac{d}{dt}\right|_{t=0} f(e^{it}z,w) = c\partial_y \ .
$$
Together with the fact that $0$ is a regular value of $f$, this shows 
that $f$ is transversal to $[0,\infty)$. Thus $f^{-1}([0,\infty))$ is a 
five-dimensional submanifold of $M$ whose boundary is $N$. 

Denote by $\mathrm{pr}_1$ the projection of $M$ onto the first factor and 
by $\mathrm{pr}_2$ the projection onto the second factor. We define the 
vector field $H_1$ on $M$ as the unique vector field with the 
property that $\mathrm{pr}_{1*}(H_1)$ is the Hopf vector field on $S^3$ 
and $\mathrm{pr}_{2*}(H_1)=0$. The vector field $H_2$ is defined  
similarly with the roles of $\mathrm{pr}_1$ and $\mathrm{pr}_2$ 
interchanged. For arbitrary constants $a$ and $b$ the vector field 
$aH_1 + bH_2$ preserves the canonical volume form on $M=S^3\x S^3$. 
The flow of this vector field at time $t$ maps $(z,w)$ to 
$(e^{iat}z,e^{ibt}w)$. If the difference $a-b$ is a rational multiple of 
$2\pi$, then all flow lines are closed, otherwise all flow lines are open. 

We consider now the flow line with starting point $(z,w)$. At time $t$ we 
find
$$
\langle z(t),w(t) \rangle = e^{i(a-b)t} \langle z,w \rangle \ . 
$$
If $a\not= b$ this means that every flow line emerging from a point of 
$M\setminus N$ intersects $f^{-1}([0,\infty))$ transversely during the time 
interval $[0,\frac{2\pi}{|a-b|}]$. Then an intersection of the flow line 
with $f^{-1}([0,\infty))$ occurs periodically after time intervals of 
length $\frac{2\pi}{|a-b|}$. The intersection is positive if $a>b$ and 
negative if $a<b$. 
Thus we have shown
$$
\lk((z,w),N)=\frac{a-b}{2\pi}
$$ 
if $(z,w)\not\in N$. 

For flow lines with starting point in $N$ the linking number with $N$ is 
not well-defined. 

For the asymptotic linking number of $aH_1+bH_2$ with $N$ we find
$$
\lk(aH_1+bH_2,N) = \frac{a-b}{2\pi} (\mathrm{vol}(S^3))^2 = 2(a-b)\pi^3 
\ .
$$
\end{ex}

\section{Linking numbers between vector fields and
measured foliations}\label{s:measures}

As before, we assume that $M$ is a closed oriented manifold
with $b_{1}(M)=0$, $\mu$ is a volume form on $M$, and $X$ a vector
field that is divergence-free with respect to $\mu$, i.~e.~such that
$L_{X}\mu=0$. Then the $(n-1)$-form $\eta=i_{X}\mu$ is closed, and
our assumptions imply that $\eta$ must be exact. Let $\alpha$ be a
primitive.

We generalise the discussion in section~\ref{s:sub} by replacing
the submanifold $N$ by an oriented codimension $2$ foliation
$\FF$ with a holonomy-invariant transverse measure $\nu$.
This defines a current
\begin{alignat*}{1}
C(\FF,\nu)\colon\Omega^{n-2}(M) &\longrightarrow\R \\
\omega &\longmapsto\int_{M}\omega\ww\nu \ ,
\end{alignat*}
where 
$$
\int_{M}\omega\ww\nu = \int_{T}(\int_{\FF}\omega) d\nu
$$
is defined by decomposing $\omega$ using a partition of unity
subordinate to a finite atlas of foliation charts for $(M,\FF)$,
integrating the summands over the plaques of $\FF$ in the charts, and then
integrating the result over the transversals $T$ using $\nu$.
The double integral is independent of the choices of charts and
partition of unity because $\nu$ was assumed to be
holonomy-invariant. See~\cite{S,CC} for the details of this
construction.  

The current $C(\FF,\nu)$ is closed, and is called the Ruelle--Sullivan
cycle of the invariant measure $\nu$. It will
play the role of the submanifold $N$ in section~\ref{s:sub}.
The assumption that $N$ be null-homologous over $\R$ is
then replaced by the assumption that the Ruelle--Sullivan
cycle is null-homologous: $[C(\FF,\nu)]=0\in H_{n-2}(M,\R)$.

The Ruelle--Sullivan cycle $C(\FF,\nu)$ is continuous with respect to the
$C^0$--to\-po\-lo\-gy on continuous forms. In general it is not possible 
to extend its domain of definition to $L^1$--forms. Nevertheless, it is 
possible to define $C(\FF,\nu)$ on integrable double forms $F(x,y)$ with 
the property that 
\begin{equation*} 
\int_{y\in M} \vert F(x,y) \vert_{y}\mu(y)
\end{equation*}
is a continuous form in the variable $x\in M$. Here we expand $F$ as a 
double form in the variables $x$ and $y$, and take the $x$-component multiplied 
by the norm of the $y$-component with respect to our fixed Riemannian metric.
Thus $\vert F(x,y) \vert_{y}$ is a differential form which is a product 
of $dx_i$, but whose coefficient function also depends on $y$. 
As we integrate along the leaves of $\FF$, we are only concerned with those 
summands of $F(x,y)$ whose degree in the variable $x$ equals the rank of the 
foliation, and we require $|F(x,y)|_y$ to be either zero or to induce the given 
orientation when restricted to the leaf of $\FF$ through $x$. 
By Proposition~\ref{p:fundL} the linking form $L$ satisfies this 
integrability condition. The form $C(\FF,\nu)(F(x,y))$ is obtained by 
performing the integrations in the Ruelle--Sullivan cycle with respect 
to the first variable. The result is an integrable form as can be shown 
with Fubini's theorem  
\begin{equation*}
\int_{y\in M} \big|C(\FF,\nu(x))\big(F(x,y)\big)\big| \ \mu(y) \le 
C(\FF,\nu(x))\left(\int_{y\in M} |F(x,y)|_{y} \ \mu(y)\right) \ .
\end{equation*} 

We define a Hopf-type integral for $X$ and $(\FF,\mu)$ by
setting
$$
H(X,\FF,\nu)=\int_{M}\alpha\ww\nu=C(\FF,\nu)(\alpha) \ .
$$
This is independent of the choice made for $\alpha$ because
the Ruelle--Sullivan cycle is assumed to be null-homologous.

We want to interpret this integral as an average of asymptotic
linking numbers of flow lines of $X$ with $\FF$. To do so, we
need again a suitable system of short paths. We use the same
notation for the paths and closed-up flow lines as before.

\begin{defn}\label{d:sps}
A {\it system of short paths} in $M$ is a set $\Sigma$ of piecewise
differentiable paths with the following properties:
\begin{enumerate}
\item For every pair of points $p,q\in M$ there is exactly one oriented
path $\sigma(p,q)\in\Sigma$ having starting point $p$ and end point $q$.
\item The paths depend continuously on their starting and
end points almost everywhere.
\item The limit
\begin{equation}\label{spscond}
\lim_{t\to\infty} \frac{1}{t} \int_{y\in\sigma(\phi_t(x),x)}
C(\FF,\nu(p))(L(p,y))= 0
\end{equation}
exists in the $L^1$--sense.
\end{enumerate}
\end{defn}
In the case when the holonomy-invariant measure $\nu$ is given by a
smooth differential form $\beta$, the Ruelle--Sullivan cycle is given by
$$
C(\FF,\nu)(\omega) = \int_M \omega\ww\beta\ .
$$
The proof of Theorem 4 in~\cite{V} generalises verbatim to this situation
and 
gives:
\begin{thm}\label{t:smoothmeasures}
If the transverse measure is given by a smooth 
ho\-lo\-no\-my-in\-va\-ri\-ant
$2$--form $\beta$, then a set of length-minimizing geodesics is a system of
short paths. 
\end{thm}
More generally, we have:
\begin{thm}\label{t:arbitrarymeasures}
Let $\FF$ be an oriented foliation with an arbitrary 
ho\-lo\-nomy-in\-va\-riant transverse
measure $\nu$. Then there exists a system of short paths.
\end{thm}
\begin{proof}
We want to generalise the proof of \thmref{t:spsN}. We use the notation 
introduced there, except that we now define $C(y)=C(\FF,\nu(p))\big(L(p,y)\big)$. 
If we can satisfy the two conditions for the choice of the base 
points $u_i\in U_i$ (with the generalised definition of $C(y)$) in the proof 
of \thmref{t:spsN}, then we can construct a system of short paths $\Sigma$ 
just as in the proof of Theorem~\ref{t:spsN}.

Fix arbitrary points $\tilde{u}_i\in U_i$ for all $1\le i \le r$ and paths 
$\tilde{\gamma}_{ij}$ without self intersection joining $\tilde{u}_i$ and 
$\tilde{u}_j$. For each pair $i\not= j$, extend the velocity vector field 
along $\tilde{\gamma}_{ij}$ to a vector field on $M$ whose time-one-flow 
transports a small ball contained in $U_i$ around $\tilde{u}_i$ to another 
ball contained in $U_j$ around $\tilde{u}_j$. Because $C(y)$ is an integrable 
form on $M$, for $\mu$--almost every starting point in a ball around 
$\tilde{u}_i$ the integral in \eqref{e:intijalt} exists. Hence for almost 
every choice of $u_i$ in the ball around $\tilde{u}_i$ we meet the first 
condition. The $\gamma_{ij}$ are the time-one-flowlines.

Also the second condition for the $u_i$ is satisfied outside of a set of 
measure zero. To see this, recall that the Ruelle--Sullivan cycle can be 
represented as a finite sum such that every summand is an integral with 
respect to a product measure. Apply Fubini's theorem and the triangle 
inequality to the integral  
\begin{multline}\label{e:rcalt}
\int_{u\in U_i} \left( \int_{\{x | n(x) = i\}} \int_{y\in \sigma(x,u)}
\left|C(\FF,\nu(p))\big(L(p,y)\big)\right|\ \mu(x)\right)\mu(u) 
\\
\le C(\FF,\nu(p))\left(\int_{s=0}^1 \int_{u\in U_i}\int_{\{x | n(x) = i\}} 
 \left|i_{\frac{\partial\sigma}{\partial s}}L\big(p,\sigma(x,u)(s)\big)\right| \mu(x)\ \mu(u) ds \right)~.
\end{multline}
If the last expression is well-defined, then by Fubini's theorem 
the expression obtained by dropping the integration with respect to $u$ is 
well-defined for $\mu$--almost every choice of $u\in U_i$. Thus we have to 
show that the form we apply the Ruelle--Sullivan cycle to is continuous. 
Consider 
\begin{equation}\label{e:int2alt}
\int_{y\in U_i}\int_{\{x | n(x) = i\}} 
\left|i_{\frac{\partial\sigma}{\partial s}}L\big(p,\sigma(x,y)(s)\big)\right| \ \mu(x)\ \mu(y) 
\end{equation}
with $p\in M$. For fixed $s\in [0,1]$ the integrand has at most a pole of order 
$n-1$ along an $n$-dimensional submanifold (the solutions of $\sigma(x,y)(s)=p$) 
in the $2n$-dimensional product manifold $U_i \x  \{x | n(x) = i\}$. This shows 
that the integral~\eqref{e:int2alt} is well-defined. Furthermore, the integral 
depends continuously on $s$ and $p$. In particular, it does exist for the 
boundary values $s=0$ and $s=1$. This shows that we apply the Ruelle--Sullivan 
cycle to a continuous form and thereby justifies the application of the 
Ruelle--Sullivan cycle in~\eqref{e:rcalt}.

So, if we choose $(u_1,\ldots,u_r)\in U_1\x\ldots\x U_r$ outside of a set of 
measure zero, we obtain a system of short paths as in the proof of 
Theorem~\ref{t:spsN}.
\end{proof}

\begin{rem}\label{r:sing}
So far we have only considered non--singular foliations. In the next section,
we will also want to use singular foliations. We therefore point out
that, as in~\cite{V}, Theorem~\ref{t:smoothmeasures} applies equally well to
singular foliations with a holonomy-invariant smooth differential form.

Theorem~\ref{t:arbitrarymeasures} can sometimes be applied to singular
foliations. For example, this can be done if the support of the transverse
measure has a neighbourhood to which the foliation extends in a nonsingular
way. In the situation of the previous section, a single
null-homologous submanifold $N\subset M$ can always be extended to a
smooth foliation of a whole neighbourhood of $N$, and we can take the
measure given by $N$, with support $N$. Then the above theorem can be
used instead of Theorem~\ref{t:spsN}.
\end{rem}

Now the linking number of $\gamma(x,t)$ with $\FF$ is defined by
generalising~\eqref{eq:link} as follows:
\begin{defn}
The linking number $\lk(\gamma,\FF,\nu)$ of a (null-homologous)
closed loop $\gamma$ in $M$ with the measured foliation $(\FF,\nu)$
is the evaluation of the Ruelle--Sullivan cycle $C(\FF,\nu)$ on a
$(n-2)$-form $\alpha$ with the property that $d\alpha$ is
Poincar\'e dual to $\gamma$.
\end{defn}
As the Ruelle--Sullivan cycle is assumed to be null-homologous,
this evaluation is independent of the choice of $\alpha$.

The following is the adaption of Proposition~\ref{p:asymptotic}
to this situation:
\begin{prop}\label{p:asymptotic2}
Let $\alpha(x,t)$ be $(n-2)$-forms with the property that
$d\alpha(x,t)$ are Poincar\'e duals for $\gamma(x,t)$. Then
the limit 
$$
\lk(x,\FF,\nu)=\lim_{t\to\infty}\frac{1}{t}\lk(\gamma(x,t),\FF,\nu)
=\lim_{t\to\infty}\frac{1}{t}C(\FF,\nu)(\alpha(x,t))
$$
exists in the $L^{1}$-sense. It is an integrable function on $M$
which does not depend on the chosen system of short paths.
\end{prop}
\begin{proof}
By the definition of $\lk(x,\FF,\nu)$ and Proposition~\ref{p:fundL} we find
\begin{align*}
\lk(x,\FF,\nu)= & \lim_{t\to\infty}\frac{1}{t}\lk(\gamma(x,t),\FF,\nu) =
               \lim_{t\to\infty}\frac{1}{t}C(\FF,\nu)(\alpha(x,t)) \\
              = & \lim_{t\to\infty}\frac{1}{t}C(\FF,\nu(p))\left(
	       \int_{y\in M} L(p,y)\ww d\alpha(x,t)(y) \right)\ .
\end{align*}
The harmonic and exact terms in \eqref{e:fundL} do not contribute because 
the Ruelle--Sullivan cycle is assumed to be null-homologous. If $p$ does 
not lie on $\gamma(x,t)$, and hence for almost every $x\in M$, we have
\begin{align*}
\lk(x,\FF,\nu)= & \lim_{t\to\infty}\frac{1}{t}C(\FF,\nu(p))\left(
              \int_{y\in\gamma(x,t)} L(p,y)\right) \\
              = & \lim_{t\to\infty}\frac{1}{t}C(\FF,\nu(p))\left(
              \int_{y\in\phi(x,t)} L(p,y)\right) \\
              = & \lim_{t\to\infty}\frac{1}{t}C(\FF,\nu(p))\left(
	      \int_0^t i_XL(p,\phi_s(x)) ds \right)\ .
\end{align*}
The second equality is true because of the definition of the system of short
paths. 
The flow of $X$ preserves the volume form $\mu$ and we can apply
the mean ergodic theorem. Hence, the limit on the right hand side of
\begin{align*}
\lk(x,\FF,\nu) = & \lim_{t\to\infty}\frac{1}{t} \int_0^t
               C(\FF,\nu(p))\big(i_XL(p,\phi_s(x))\big) ds
\end{align*} 
exists in the $L^1-$sense and represents an integrable function on $M$. It
does not depend on the system of short paths.
\end{proof} 

Using this, we can finally define the average asymptotic linking
number of the vector field $X$ with the measured foliation $(\FF,\nu)$
by setting
$$
\lk(X,\FF,\nu)=\int_{M}\lk(x,\FF,\nu)\mu \ .
$$
Theorem~\ref{t:sub} generalises as follows:
\begin{thm}\label{t:linkmeasure}
Let $M$ be a closed oriented $n$-manifold with $b_1(M)=0$.
Let $X$ be a divergence-free vector field on $M$, and $\FF$ an
oriented codimension $2$ foliation with a transverse measure $\nu$
whose Ruelle--Sullivan cycle is null-homologous. Then the average
asymptotic linking number of the orbits of $X$ with $(\FF,\nu)$
exists and equals a Hopf-type integral:
$$
\lk(X,\FF,\nu)=H(X,\FF,\nu) \ .
$$
\end{thm}
\begin{proof}
We use the calculations in the proof of Proposition~\ref{p:asymptotic2} and
the mean ergodic theorem. Since $b_1(M)=0$, $b_{n-1}(M)$ also vanishes.
The $(n-1)$--form $i_X\mu$ is closed and hence exact. Choose an
$(n-2)$--form
$\alpha_X$ such that $d\alpha_X=i_X\mu$. By the mean ergodic theorem and
\propref{p:asymptotic2} 
\begin{align*}
\lk(X,\FF,\nu) = & \int_{x\in M}\lk(x,\FF,\nu)\mu(x) \\
               = & \int_{x\in M}C(\FF,\nu(p))(i_XL(p,x))\mu(x)\ .
\end{align*}       
Locally, the Ruelle--Sullivan cycle is given by a product measure, hence we
can apply Fubini's theorem. We find
\begin{align*}
\lk(X,\FF,\nu) = & C(\FF,\nu(p))\left( \int_{x\in M}i_XL(p,x)
               \ww\mu(x)\right) \\  
               = & C(\FF,\nu(p))\left(\int_{x\in M}L(p,x)\ww i_X\mu(x)
               \right) \\
               = & C(\FF,\nu(p))(\alpha_{X}) =  H(X,\FF,\nu) \ ,
\end{align*}
where the penultimate equality is due to Proposition~\ref{p:fundL} and 
Stokes's theorem. 
\end{proof}

\begin{rem}
The discussion in this section reduces to that of the previous
section in the case that the invariant measure $\nu$ is given by a
closed leaf $N$. The rest of the foliation $\FF$ then plays no role.
\end{rem}

\begin{rem}\label{r:Khesin}
For a smooth foliation with a holonomy-invariant measure given by a 
smooth exact $2$-form, one can prove Theorem~\ref{t:linkmeasure} 
using Arnold's definition of a system of short paths and the Birkhoff 
ergodic theorem as in~\cite{A}, rather than the mean ergodic theorem as 
above. This was done by Khesin in~\cite{Kh}.
\end{rem}

Here are two examples for linking numbers between measured foliations
and divergence-free vector fields.

\begin{ex}
Let $\FF$ be the Reeb foliation on $S^3$. It has exactly one closed leaf, which is 
diffeomorphic to $T^2$. All other leaves are diffeomorphic to $\R^2$ and have 
linear growth. By a construction going back to Plante, cf.~\cite{CC,S}, any leaf 
$\mathcal{L}$ of subexponential growth defines an invariant measure 
$\mu_{\mathcal{L}}$ with support contained on a union of minimal sets 
in $\overline{\mathcal{L}}$.
For the Reeb foliation one can show easily that, up to 
a factor of $2\pi$, the measure $\nu_{\mathcal{L}}$ defined by an 
open leaf equals that defined by the unique closed leaf $T^{2}$.

Consider now the product of two Reeb foliations on $S^3\x S^3$. This foliation 
$\FF\x\FF$ has codimension two and contains exactly one closed leaf diffeomorphic 
to $T^4$. If $\mathcal{L}_{1}$ and $\mathcal{L}_{2}$ are open leaves of 
$\FF$, then the leaf $\mathcal{L}_{1}\x\mathcal{L}_{2}$ in the product foliation 
defines the holonomy-invariant transverse measure 
$\nu_{\mathcal{L}_{1}}\times\nu_{\mathcal{L}_{2}}$ for $\FF\times\FF$. 
By the discussion above we have the equality
\begin{equation*}
\nu_{\mathcal{L}_{1}}\times\nu_{\mathcal{L}_{2}} = \frac{1}{4\pi^2}\nu_{T^4} \ .
\end{equation*}
Thus the linking number of a divergence-free vector field with 
$(\FF\x\FF, \nu_{\mathcal{L}_{1}}\times\nu_{\mathcal{L}_{2}})$ is, up to 
a factor of $4\pi^{2}$, the same as that with the submanifold 
$T^4\subset S^{3}\times S^{3}$ given by the closed leaf.
\end{ex} 

\begin{ex}
Next we consider a singular foliation on $\cp$. In order to obtain a
null-homologous Ruelle--Sullivan cycle we will use signed measures, rather 
than measures. Consider the affine embedding
\begin{alignat*}{1}
\C^2 & \hookrightarrow \cp \\
(z_0,z_1) & \mapsto [1:z_0:z_1] \ .
\end{alignat*}
We foliate $\C^2$ by the product foliation $\C\x\C$ and consider the leaves 
$\mathcal{L}_1=\{1\}\x\C$ and $\mathcal{L}_{-1}=\{-1\}\x\C$. 
Although these leaves are not closed, the counting measure is a well-defined 
holonomy-invariant transverse measure since the closure of $\mathcal{L}_1$,
respectively $\mathcal{L}_{-1}$, is the union of the leaf with $\{[0:0:1]\}$. 
From now on we equip $\mathcal{L}_1$ with the counting measure and 
$\mathcal{L}_{-1}$ with the negative counting measure, i.~e.~the signed measure 
which counts each intersection point with $-1$. We denote this signed measure 
by $\delta_1 - \delta_{-1}$. 

The corresponding Ruelle--Sullivan cycle is 
\begin{alignat*}{1}
C(\FF, \delta_1 - \delta_{-1}) : \Omega^2(M) & \lra \R \\
\omega      & \longmapsto \int_{\mathcal{L}_1} \omega - 
\int_{\mathcal{L}_{-1}} \omega \ .
\end{alignat*}
Here $C(\FF, \delta_1 - \delta_{-1})$ is null-homologous by construction.

Consider $N=\left\{(z,w)\in\C^2\subset\cp | z\in[-1,1]\subset\R \right\}$. 
The boundary of this submanifold of $\C^2$ is exactly the union of the two 
leaves $\mathcal{L}_1$ and $\mathcal{L}_{-1}$, the latter with the reversed 
orientation. The corresponding current 
\begin{alignat*}{1}
\Omega^3(\cp) & \lra \R \\
\omega        & \longmapsto \int_N \omega
\end{alignat*}
has the boundary $C(\FF, \delta_1 - \delta_{-1})$.

Let $\mu$ be a volume form on $\cp$ and $X$ a divergence--free  vector field. 
Let $\al$ be a primitive of $i_X\mu$. Then the linking number of $X$ with 
$(\FF, \delta_1 - \delta_{-1})$ is 
\begin{align*}
\mathrm{lk}(X,\FF,\delta_1 - \delta_{-1}) & = \int_{\mathcal{L}_1} \al - \int_{\mathcal{L}_{-1}} \al \\
                                          & = \int_N i_X\mu~.
\end{align*}
This is exactly the flux of $X$ through $N$.
\end{ex}

\section{Godbillon-Vey invariants as linking numbers}\label{s:GV}

Recall that the tangent distribution of a smooth codimension $1$
foliation $\FF$ with trivial normal bundle on a manifold $M$ is
the kernel of a non-vanishing $1$-form $\alpha$ which, by the
Frobenius theorem, satisfies
\begin{equation}\label{eq:da}
d\alpha = \alpha \ww \beta
\end{equation}
for some $1$-form $\beta$. The $3$-form $\beta\ww d\beta$
is closed, and its cohomology class $GV(\FF)\in H^{3}(M,\R)$ is
independent of the choices made for $\alpha$ and $\beta$. This
is the Godbillon-Vey invariant~\cite{GV} of $\FF$.

In the case that $M$ is closed, oriented and $3$-dimensional,
$GV(\FF)$ is equivalent to the Hopf integral
\begin{equation}\label{eq:H}
\int_{M}\beta\ww d\beta \ .
\end{equation}
If we choose an arbitrary volume form $\mu$ on $M$ and define
a vector field $X$ in $M$ by the formula $i_{X}\mu = d\beta$,
then $X$ is divergence-free with respect to $\mu$, and~\eqref{eq:H}
is just the Hopf invariant $H(X,X)$. If $M$ is a $\R$-homology
sphere, then by the ``helicity theorem'' due to Arnold~\cite{A} and
the second author~\cite{V}, this Hopf invariant can be interpreted
as the average asymptotic self-linking number of the orbits of $X$.

Assume now that we have a smooth $1$-parameter family of
smooth codimension $1$ foliations $\FF_{t}$ with trivial normal
bundles. Then~\eqref{eq:da} still holds, but now $\alpha$ and
$\beta$ are functions of the parameter $t\in\R$. We denote the
time derivatives by a dot. It was shown by the first author~\cite{K}
that for every $t$ the $4$-form $(\dot\beta\ww\beta\ww d\beta)(t)$ is
closed, and that its cohomology class $TGV(\FF_{t})\in H^{4}(M,\R)$
is a well-defined invariant of the family $\FF_{t}$ that is
independent of choices.

In the case that $M$ is closed, oriented and $4$-dimensional,
$TGV(\FF_{t})$ is equivalent to the integral
\begin{equation}\label{eq:TGV}
\int_{M} \dot\beta\ww\beta\ww d\beta(t) \ .
\end{equation}
We want to give an interpretation of this as an average asymptotic
linking number of a suitable vector field and a measured codimension
$2$ foliation.

Choose an arbitrary volume form $\mu$ on $M$ and define a 
time-dependent vector field $X$ by the formula 
$i_{X}\mu = d(\dot\beta\ww\beta)$.

Differentiating~\eqref{eq:da}, we see that $d\beta\ww d\beta =0$
because $d\beta$ is decomposable. Thus, on the open set in $M$ where
$d\beta$ does not vanish, its kernel is a $2$-dimensional distribution
which is integrable because the defining form is closed. We denote by
$\GG$ the singular codimension $2$ foliation defined by $d\beta$. Note
that the exact form $d\beta$ defines an invariant transversal measure
for $\GG$ whose Ruelle--Sullivan cycle is null-homologous. Now we
recognise~\eqref{eq:TGV} as the Hopf-type integral associated to
the vector field $X$ (at time $t$) and the measured foliation 
$(\GG,d\beta)$ (at the same time $t$).

To interpret this as an average asymptotic linking number
according to Theorem~\ref{t:linkmeasure}, we need to assume
$b_{1}(M)=0$. As $M$ is closed, oriented and $4$-dimensional,
it follows that the Euler characteristic of $M$ is positive, and
so there cannot be any non-singular codimension $1$ foliation $\FF$
on $M$. However, as we have ended up with only a singular foliation
for $\GG$, there is no harm in allowing $\FF$ to be singular as
well. So we just assume that $\FF_{t}$ is the kernel
of a time-dependent $1$-form $\alpha$ satisfying~\eqref{eq:da}, but allow
$\alpha$ to have zeroes. Under certain technical assumptions on 
$\alpha$, cf.~the appendix, the definition of $TGV(\FF_{t})$ goes through
as in the non-singular case, and if $b_{1}(M)=0$, then according
to Theorem~\ref{t:linkmeasure} the integral of $TGV(\FF_{t})$ over
$M$ is the average asymptotic linking number $\lk(X,\GG,d\beta)$, with
$X$ and $\GG$ defined as above. We can apply
Theorem~\ref{t:linkmeasure} to the singular foliation $\GG$ because
the holonomy-invariant measure is given by a smooth form, see
Remark~\ref{r:sing}.

If $\FF_{t}$ is a $1$-parameter family of codimension $q$ foliations
on $M$ defined by a decomposable $q$-from $\alpha$, then setting $d\alpha
=\alpha\wedge\beta$ as above, we can consider
$\dot\beta\wedge\beta\wedge (d\beta)^{q}$. It was proved in~\cite{K}
that this is closed, and that its cohomology class $TGV(\FF_{t})\in
H^{2q+2}(M,\R)$ is a well-defined invariant of the family $\FF_{t}$.
One might be tempted to think that when $M$ is closed oriented of
dimension $2q+2$, then this invariant should always be an average asymptotic
linking number, as was proved above for the case $q=1$. However, this
does not seem to be the case. We always have $(d\beta)^{q+1}=0$,
showing that $d\beta$ has at least a $2$-dimensional kernel. If
$TGV(\FF_{t})$ doesn't vanish, then $(d\beta)^{q}\neq 0$ on an open
set in $M$, so that the kernel of $d\beta$ is exactly of rank $2$ (on
this open set). This means that whenever $q>1$, the codimension of
the cokernel of $d\beta$ is strictly larger than $2$, and so there
can be no linking number with $1$-dimensional flow lines.

\section*{Appendix:\\ Singular foliations and forms with the division property}

In the last section we have had to consider codimension one foliations 
defined by one-forms $\alpha$ which have zeroes. In this appendix we 
explain how the definitions of the Godbillon-Vey invariant $GV$ and of the 
invariant $TGV$ of families extend to this situation. As in the rest 
of this paper, all forms are smooth of class $C^{\infty}$.

Following Moussu~\cite{M}, we consider the following:
\begin{defn}
    A one-form $\alpha$ has the division property 
    if $\omega\wedge\alpha=0$ for a smooth $\omega\in\Omega^{k}(M)$ 
    with $0 < k < \dim (M)$ implies $\omega=\alpha\wedge\beta$ for 
    some smooth $\beta\in\Omega^{k-1}(M)$.
    \end{defn}
Note that the division property 
then implies 
that $\beta$ is unique up to the addition of multiples of $\alpha$.

Nowhere vanishing $1$-forms have the division property,
as do some classes of singular $1$-forms. 
For example, it is easy to see that if near every point where it 
vanishes, $\alpha$ is locally the differential of a Morse function, 
then it has the division property. More generally, Moussu~\cite{M} 
proved that $1$-forms with only algebraically isolated zeros have the 
division property. 
Algebraic isolation of the zeros means that in local coordinates at a 
zero the coefficient functions of the one-form span an ideal of finite 
codimension in the algebra of germs of functions.

The calculation showing the well-definedness of the Godbillon-Vey 
invariant~\cite{GV} of a codimension one foliation does not require 
the defining form $\alpha$ to be non-zero, but rather requires only 
that it have the division property. 
Thus we have:
\begin{thm}
    Let $\alpha$ be a smooth $1$-form with $\alpha\wedge d\alpha=0$, and 
    $\FF$ its kernel foliation. If $\alpha$ has the division property 
    then there is a $1$-form $\beta$ with
\begin{equation}\label{eq:daA}
d\alpha = \alpha \ww \beta \ .
\end{equation}
The $3$-form $\beta\ww d\beta$ is closed and its cohomology class 
$GV(\FF)\in H^{3}(M,\R)$ is independent of the choices 
made for $\alpha$ and $\beta$, as long as $\alpha$ is changed only by 
multiplication with a nowhere vanishing function.
\end{thm}
    
This generalises to families in the following way:
\begin{thm}
Let $\alpha$ be a smooth $1$-form on $M$ with $\alpha\wedge d\alpha=0$ and 
which depends smoothly on a parameter $t\in\R$. Let $\FF_{t}$ be the family 
of kernel foliations. If for every $t$ the form $\alpha$ has the division 
property,
then the $4$-forms
$(\dot\beta\wedge\beta\ww d\beta)(t)$ are closed, and their cohomology 
classes $TGV(\FF_{t})\in H^{4}(M,\R)$ are independent of the choices made for 
$\alpha$ and $\beta$, as long as $\alpha$ is changed only by 
multiplication with a nowhere vanishing function.
\end{thm}
\begin{proof}
    Once we have the forms $\beta$ and $\dot\beta$, checking the 
    claims in the theorem is done exactly as in the nonsingular case, 
    compare~\cite{K}.
    
    In the case when $\alpha$ has no zeros, it is easy to see that if 
    $\alpha$ depends smoothly on $t$, then one can also choose $\beta$ 
    in~\eqref{eq:daA} to depend smoothly on $t$. In the singular case, the 
    smooth dependence of $\beta$ on $t$ is less obvious. However, if 
    $\alpha(0)$ has only algebraically isolated zeros, then the same is 
    true for $\alpha(t)$ for all $t$ sufficiently close to $0$, and in this 
    case $\beta$ can be chosen to depend smoothly on $t$ by adapting 
    Moussu's argument~\cite{M}.
    
    One can also give an alternative treatment for all divisible forms 
    as follows.
    Recall from~\cite{K} the identity
\begin{equation}\label{eq:alt}
    \dot\bb\ww\bb\ww d\bb= d(\dot\al\ww\bb\ww\cc )
-\dot\al\ww\bb\ww\al\ww\delta \ ,
\end{equation}
where $\cc$ and $\delta$ are $1$-forms constructed by successive 
exterior differentiation of~\eqref{eq:daA} and the division property 
of $\alpha$. This shows that $TGV(\FF_{t})$ 
can be defined using only the smooth dependence of $\alpha$ on $t$, 
and well-definedness follows as in~\cite{K} using the division property.
    \end{proof}


\bibliographystyle{amsplain}

\end{document}